\documentstyle[12pt]{article}
\renewcommand{\theequation}{\thesection.\arabic{equation}}
\newcounter{subequation}[equation]
\makeatletter

\expandafter\let\expandafter
\reset@font\csname reset@font\endcsname

\def\subeqnarray{\arraycolsep1pt
    \def\@eqnnum\stepcounter##1{\stepcounter{subequation}%
        {\reset@font\rm(\theequation\alph{subequation})}}
\jot5mm     \eqnarray}

\newfont{\blackb}{msbm10 scaled\magstep1}
\def\Bbb#1{\hbox{\blackb #1}}
\def \a{\alpha}
\def \b{\beta}
\def \g{\gamma}
\def \d{\delta}
\def \e{\varepsilon}
\def \l{\lambda}
\def \L{\Lambda}
\def \k{\kappa}
\def \t{\triangle}

\def \p{\ \>}
\def \G{\Gamma}
\def \big{\bigtriangledown}
\makeatother
\newcommand{\newsection}[1]{
\pagebreak[3]
\addtocounter{section}{1}
\setcounter{equation}{0}
\setcounter{subsection}{0}
\setcounter{footnote}{0}

\begin{flushleft}
{\Large\bf \thesection. #1}
\end{flushleft}
\nopagebreak
\medskip
\nopagebreak}

\begin{document}
\baselineskip 18pt
\vspace{1cm}
\begin{center}
\begin{Large}
{\bf Dirac operator on $\k$-Minkowski space, bicovariant differential
calculus and deformed $U(1)$ gauge theory}
\end{Large}
\\
\vspace{1cm}
{\bf P.N. Bibikov}
\vskip10pt
{\it St.Petersburg Branch of Steklov Mathematical
Institute}\vskip0pt
{\it Fontanka 27, St.Petersburg, 191011, Russia}
\\
\vskip10pt
\end{center}

\begin{abstract}                
Derivation of $\k$-Poincare bicovariant commutation relations between
coordinates and 1-forms on $\k$-Minkowski space is given using Dirac operator
and Allain Connes formula. The deformed $U(1)$ gauge
theory and appearance of an additional spin 0 gauge field is discussed.
\end{abstract}
\setcounter{section}{-1}

\newpage

\newsection{Introduction}

Noncommutative geometry suggested by A. Connes \cite{1} attracts these days
a great interest of many researchers. Besides a lot of mathematical
applications \cite{1},\cite{13} it also may be considered as a natural
framework for quantisation of space and time \cite{11}.
One of the most promising results in this direction is
the approach to gauge field theory developed in \cite{2} where the Standard
Model of gauge interaction was obtained from noncommutativity of space-time.

The basic notion of the approach studied in \cite{1}, \cite{2} is the Connes
triple (${\cal A}, {\cal H}, D$) where ${\cal A}$ is in general
framework a noncommutative $^*$-algebra which is considered as algebra of
operators in the Hilbert space ${\cal H}$. $D$ is
a linear possibly unbounded operator in ${\cal H}$ with $D^*=-D$.

In classical case when ${\cal A}={\rm Fun}(M)$ is commutative algebra of
functions on the differential manifold $M$
\begin{equation}
D=\g^{\mu}\partial_{\mu}\,,
\end{equation}
is the usual Dirac operator.

In (0.1) $\partial_i$ are local derivatives and $\g^i=\g^i(x)$ are generators
of local Clifford algebres
\begin{equation}
\g^{\mu}(x)\g^{\nu}(x)+\g^{\nu}(x)\g^{\mu}(x)=2g^{\mu\nu}(x)\,,
\end{equation}
where $g^{\mu\nu}(x)$ are local components of metric tensor.

Noncommutative differential calculus on the Fun($M$) is defined by
introduction exterior derivative operator which we shall denote by $d_c$.
For each $f\in{\rm Fun}(M)$ it has the form
\begin{equation}
d_cf=[D,f]\,,
\end{equation}

According to (0.1) formula (0.3) gives the following result
\begin{equation}
d_cf=\partial_{\mu}(f)\g^{\mu}
\end{equation}

Correspondence of the definition (0.3) with the usual external derivative
\begin{equation}
df=\partial_i(f)dx^i
\end{equation}
follows from the isomorhism between $\G$ and $Cl({\rm Fun}(M))$ where $\G$ is
the space of differential 1-forms over $M$ or the space of sections of
cotangent bundle $T^*(M)$ and $Cl({\rm Fun}(M))$ is the space of sections
of $Cl(M)$ the Clifford bundle over $M$ \cite{4}. This may be expressed by
commutative diagram
\begin{eqnarray}
\begin{array}{ccc}
{\rm Fun}(M)&\stackrel{df=\partial_i(f)dx^i}\longrightarrow &\Gamma\\
\downarrow\lefteqn{id}&&\downarrow\lefteqn{dx^i\rightarrow\g^i}\\
{\rm Fun}(M)&\stackrel{d_cf=[D,f]}\longrightarrow & Cl({\rm Fun}(M))
\end{array}
\end{eqnarray}

The fiber of $Cl(M)$ corresponding to $x\in M$ is a
Clifford algebra defined by generators $\g^{\mu}(x)$ and relations (0.2).
$Cl({\rm Fun}(M))$ is a bimodule over ${\rm Fun}(M)$ generated by all sums of
the form
\begin{equation}
\sum_{i}f_i[D,g_i]\,.
\end{equation}

According to isomorphism (0.6) gauge connection 1-form $A_{\mu}dx^{\mu}$ which
is used in construction of pure gauge action and the gauge interaction term
$A_{\mu}\g^{\mu}$ in Dirac equation for spinor field have similar geometrical
intrpretations. So studying deformations of field theory in quantum spaces it
is natural to suppose that the diagram (0.6) has analog also in noncommutative
case. We shall write corresponding diagram in the form,
\begin{eqnarray}
\begin{array}{ccc}
{\cal A}&\stackrel{d:f\rightarrow df}\longrightarrow & \Gamma\\
\downarrow\lefteqn{id}&&\downarrow\lefteqn{df\rightarrow d_cf}\\
{\cal A}&\stackrel{d_cf=[D,f]}\longrightarrow & Cl({\cal A})
\end{array}
\end{eqnarray}
where $Cl({\cal A})$ is bimodule over ${\cal A}$ generated by all sums of
the form (0.7).

For the most interesting examples of noncommutative manifolds
studied in the context of quantum group theory operator $d$ and corresponding
spaces of quantum 1-forms are defined axiomatically \cite{3}. So in this
case the condition of commutativity diagram (0.8) gives strong restrictions on
the possible form of Dirac operator.

In this framework the Dirac equation for a massless spinor field coupled
with gauge potential has the form
\begin{equation}
(D+gV)\psi=0\,,
\end{equation}
where $\psi\in{\cal H}$ and $V$ is a noncommutative analog of
$iA_{\mu}\g^{\mu}$ and $g$ is a gauge charge. According
to isomorphism between quntum Clifford and quantum cotangent bundles supposed
by (0.8) it correspons to $\omega$ the gauge connection quantum 1-form which
is the noncommutative analog of $iA_{\mu}dx^{\mu}$.

We suggest deformed $U(1)$ gauge transformation of spinorial fields in the form
\begin{equation}
\psi\rightarrow U\psi\,,
\end{equation}
where $U$ is a unitary element of ${\cal A}$
\begin{equation}
UU^*=U^*U=1\,.
\end{equation}
(Additional restriction on $U$ will be discussed in the last section).
The transformation (0.10) for $\psi$ is compatible with the following
transformation for $V$
\begin{equation}
V\rightarrow UVU^*+U[D,U^*]
\end{equation}
that according to (0.8) is equivalent to the following law for $\omega$
\begin{equation}
\omega\rightarrow\tilde\omega=U\omega U^*+UdU^*\,.
\end{equation}

In the present paper we construct Dirac operators for quantum spaces
admitting quantum group coaction. Bicovariant differential calculus on
these spaces can be constructed according to \cite{3}. Commutativity of the
diagram (0.8) follows from general theory of quantum groups. As an example
we consider the case of $\k$-Minkowski spase ${\cal M}_{\k}$ which 
corresponds to the $\k$-Poincare group and was intensively  studied in the 
last few years. Description of different Minkowski space deformations
is given in \cite{11}. (Also Dirac operators on the quantum $SU(2)$ group
and the quantum sphere were obtained in \cite{8} and \cite{12}).   
In paper \cite{10} a Dirac operator on ${\cal M}_{\k}$ was defined. In this
paper we present the one parametric family of Dirac operators on
${\cal M}_{\k}$ and among them lies the one suggested in \cite{10}. 

The paper is organised as follows. 
In sect. 1 we study according to \cite{5},\cite{6} and \cite{7} 
differential geometry on ${\cal M}_{\k}$.
Invariant Klein-Gordon operator on ${\cal M}_{\k}$ is constructed as a
bilinear combnation of generators of covariant algebra of vector fields,
expressed from standard generators of $\k$-Poincare quantum algebra.
Construction of exterior differential needs introduction of quantum derivatives
which also are elements of quantum Poincare algebra.
In sect. 2 we construct on ${\cal M}_{\k}$ Dirac operator and show 
commutativity of (0.8).
In sect. 3 we suggest on ${\cal M}_{\k}$ equations of deformed electrodynamics.
We also discuss deformation of corresponding gauge invariant action 
and illustrate how the existence of extra dimension of space of
quantum differential 1-forms can explain appearance of spin 0 gauge field.

Everywhere in the paper we use Einstein rules of summation. In first two
sections greece indices $\a,\mu,\nu$ numerate space-time components and
take value 0,1,2,3 however latin indices $m,n,$ numerate only the space
components and take value 1,2,3. Everywhere $g^{\mu\nu}$ means the
Minkowski space metric tensor ($1,-1,-1,-1$).

The author is very greateful for P. P. Kulish for posing the problem 
and useful informative discussions.

\newsection{$\k$-Poincare group and $\k$-Minkowski space}

The $\k$-Poincare quantum group was introduced in \cite{6} (see also \cite{5},
\cite{7}) and in one of
equivalent forms it represents as a $^{*}$-Hopf algebra generated by
hermitian elements $\L_{\mu}^{\p\nu}$, $a^{\mu}$ and relations
\begin{eqnarray}
[a^{\mu},a^{\nu}]&=&\frac{i}{\k}(\d_{0}^{\mu}a^{\nu}-\d_0^{\nu}a^{\mu})\,,
\nonumber \\
{[}\L_{\mu}^{\p\nu},\L_{\a}^{\p\b}{]}&=&0\,,\nonumber\\
{[}\L_{\mu}^{\p\nu},a^{\a}{]}&=&\frac{i}{\k}[(\d_{0}^{\nu}-\L_{0}^{\p\nu})
\L_{\mu}^{\p\a}+g^{\nu\a}(\d_{\mu}^0-\L_{\mu}^{\p0})]\,,\nonumber\\
\t(\L_{\mu}^{\p\nu})&=&\L_{\mu}^{\p\a}\otimes\L_{\a}^{\p\nu}\,,\nonumber\\
\t(a^{\mu})&=&a^{\nu}\otimes\L_{\nu}^{\p\mu}+1\otimes a^{\mu}\,,\\
S(\L_{\mu}^{\p\nu})&=&\L^{\nu}_{\p\mu}=g_{\mu\,\a}g^{\nu\,\b}\L_{\b}^{\p\a}\,,
\nonumber\\
S(a^{\mu})&=&-a^{\nu}\L^{\mu}_{\p\nu}\,,\nonumber\\
\e(\L_{\mu}^{\p\nu})&=&\d_{\mu}^{\nu}\,,\nonumber\\
\e(a^{\mu})&=&0\,,\nonumber
\end{eqnarray}

The ${\cal P}_{\k}$ may be regarded as the quantum
symmetry group of $\k$-Minkowski space ${\cal M}_{\k}$, which is defined
by four hermitian generators $x^{\mu}$ and relations
\begin{equation}
[x^{\mu},x^{\nu}]=\frac{i}{\k}(\d_{0}^{\mu}x^{\nu}-\d_{0}^{\nu}x^{\mu})\,.
\end{equation}

The corresponding right ${\cal P}_{\k}$ coaction is
\begin{equation}
\Phi_R(x^{\mu})=x^{\nu}\otimes\L_{\nu}^{\p\mu}+1\otimes a^{\mu}
\end{equation}
(the left comodule structure also can be defined \cite{5},\cite{6}).

Coproduct, counit, and antipode \cite{7}
\begin{eqnarray}
\t(x^{\mu})&=&1\otimes x^{\mu}+x^{\mu}\otimes 1\,,\nonumber\\
\e(x^{\mu})&=&0\,,\\
S(x^{\mu})&=&-x^{\mu}\,,\nonumber
\end{eqnarray}
also define on ${\cal M}_{\k}$ structure of cocommutative Hopf algebra.
As it was shown in \cite{7} correspondence $a^{\mu}\rightarrow x^{\mu}$ and
$\L^{\mu}_{\p\nu}\rightarrow1$ defines a Hopf algebra homomorphism from 
${\cal P}_{\k}$ to ${\cal M}_{\k}$. 

Its Hopf dual ${\cal M}_{\k}^*$ is defined by four hermitian generators
$P_{\mu}$ and relations
\begin{eqnarray}
[P_{\mu},P_{\nu}]&=&0\,,\nonumber\\
\t(P_0)&=&P_0\otimes 1+1\otimes P_0\,,\nonumber\\
\t(P_m)&=&P_m\otimes 1+e^{-\frac{P_0}{\k}}\otimes P_m\,,\nonumber\\
S(P_0)&=&-P_0\,,\\
S(P_m)&=&-e^{P_0/\k}P_m\,,\nonumber\\
\e(P_{\mu})&=&0\,.\nonumber
\end{eqnarray}
and the pairing $(\cdot,\cdot):{\cal M}_{\k}^*\otimes{\cal M}_{\k}\rightarrow
{\Bbb C}$ is given by
\begin{equation}
i(P_{\mu},x^{\nu})=\d_{\mu}^{\nu}\,.
\end{equation}

As it was shown in \cite{7} ${\cal M}^*_{\k}$ is a Hopf subalgebra
of the quantum Poincare algebra which is dual to ${\cal P}_{\k}$. The 
corresponding pairing between ${\cal M}^*_{\k}$ and ${\cal P}_{\k}$ is given 
by formula
\begin{equation}
i(P_{\mu},a^{\nu})=\d_{\mu}^{\nu}\,.
\end{equation}

And according to this pairing ${\cal M}^*_{\k}$ acts on
${\cal M}_{\k}$ from the left.
\begin{eqnarray}
\pi(x)&=&((id\otimes\pi),\Phi_R(x))\,,\\
&&\pi\in{\cal M}^*_{\k}\,,\quad x\in {\cal M}_{\k}\,.\nonumber
\end{eqnarray}

Considering elements of ${\cal M}_{\k}$ as left multiplication operators we
may obtain according to (1.5) the following relations
\begin{eqnarray}
[P_0,x^{\mu}]&=&\frac{1}{i}\d_0^{\mu}\,,\nonumber\\
{[}P_m,x^0{]}&=&\frac{i}{\k}P_m\,,\\
{[}P_m,x^n{]}&=&\frac{1}{i}\d_m^n\,.\nonumber
\end{eqnarray}

Elements
\begin{eqnarray}
e^4&=&i\k({\rm ch}\frac{P_0}{\k}-\frac{1}{2\k^2}e^{P_0/\k}\vec P^2)\,,
\nonumber\\
e^0&=&i\k({\rm sh}\frac{P_0}{\k}+\frac{1}{2\k^2}e^{P_0/\k}\vec P^2)
\,,\\
e^m&=&-ie^{P_0/\k}P_m\,,\nonumber
\end{eqnarray}
satisfy following commutation relations with elements of ${\cal M}_{\k}$
\begin{eqnarray}
[e^{\mu},x^{\nu}]&=&\frac{i}{\k}(g^{0\mu}e^{\nu}-g^{\mu\nu}e^0-g^{\mu\nu}e^4)
\,,\nonumber\\
{[}e^4,x^{\mu}{]}&=&-\frac{i}{\k}e^{\mu}\,,
\end{eqnarray}
and the additional relation
\begin{equation}
\Box_{\k}\equiv e_{\mu}e^{\nu}=g^{\mu\nu}e_{\mu}e_{\nu}=
\k^2+(e^4)^2\,.
\end{equation}

Eqs. (1.11), (1.12) are invariant under the right ${\cal P}_{\k}$ coaction 
which on
${\cal M}_{\k}$ has the form (1.3) and on the elements (1.10) is defined by
\begin{equation}
\Phi_R(e_{\mu})=e_{\nu}\otimes\L^{\nu}_{\p\mu}\,,\qquad
\Phi_R(e^4)=e^4\otimes1\,.
\end{equation}

So according to the general approach \cite{9},\cite{14} we may consider the 
joint algebra generated by $x^{\mu}$, $e_{\mu}$, $e^4$ and relations
(1.2), (1.11) and (1.12) as the algebra of vector fields on ${\cal M}_{\k}$.

Element $\Box_{\k}$ from (1.12) is invariant under (1.13) and playes a role
of massless Klein-Gordon operator on ${\cal M}_{\k}$ \cite{10}.

Quantum De Rham complex on ${\cal M}_{\k}$ with satisfied Leibnitz rule
was constructed in \cite{5}. The space of 1-forms $\Gamma$ is generated by
\begin{equation}
\tau^{\mu}=dx^{\mu}\,,\qquad\tau^4=\frac{i\k}{4}([\tau^{\mu},x_{\mu}]
+\frac{3i}{4}\tau^0)\,,
\end{equation}
and relations
\begin{eqnarray}
[\tau^{\mu},x^{\nu}]&=&\frac{i}{\k}(g^{0\mu}\tau^{\nu}-g^{\mu\nu} \tau^0-
g^{\mu\nu}\tau^4)\,,\nonumber\\
{[}\tau^4,x^{\mu}{]}&=&-\frac{i}{\k}\tau^{\mu}\,.
\end{eqnarray}

External algebra relations and external derivative are given by ($i,j=0,...4$)
\begin{equation}
\tau^i\wedge\tau^j=-\tau^j\wedge\tau^i\,,\qquad d\tau^i=0\,.
\end{equation}

Eqs. (1.14), (1.15) and (1.16) are invariant under the right 
${\cal P}_{\k}$-coaction
given on elements of ${\cal M}_{\k}$ by (1.3) and on elements $\Gamma$ by
\begin{equation}
\Phi_R(\tau^{\mu})=\tau^{\nu}\otimes\L_{\nu}^{\p\mu}\,,\qquad\Phi_R(\tau^4)=
\tau^4\otimes1\,.
\end{equation}
(In \cite{5} the left variant of (1.17) was presented).

It is easy to see from (1.15) that for every $a\in{\cal M}_{\k}$
\begin{equation}
s^2a=as^2\,,
\end{equation}
where the metric form $s^2\in\Gamma\otimes\Gamma$ defined by
\begin{equation}
s^2=\tau_{\mu}\otimes\tau^{\mu}-\tau^4\otimes\tau^4
\end{equation}
is invariant under the right ${\cal P}_{\k}$ coaction on $\Gamma\otimes\Gamma$

According to (1.12) and (1.18) we may define corresponding to the $e^4$ and
$\tau^4$ component of metric tensor.
\begin{equation}
g^{44}=g_{44}=-1
\end{equation}

Commutation relations between 1-forms and elements of ${\cal M}_{\k}$ also
may be represented in the standard form \cite{3} ($i,j=0,1,2,3,4$)
\begin{equation}
\tau^ia=f^i_{\p j}(a)\tau^j\,,
\end{equation}
where $\tau^i_k$ are linear operators $\tau^i_k:{\cal A}\rightarrow{\cal A}$.
From $\tau^i(ab)=(\tau^ia)b$ follows that
\begin{equation}
f^i_{\p k}(ab)=f^i_{\p j}(a)f^j_{\p k}(b)\,,
\end{equation}
In the most interesting case when all $f^i_{\p j}\in{\cal M}^*_{\k}$ so that
their action on elements of ${\cal A}$ is given by (1.8) that is
equivalent to
\begin{equation}
\t(f^i_{\p k})=f^i_{\p j}\otimes f^j_{\p k}\,.
\end{equation}
Taking
\begin{eqnarray}
f^0_{\p0}&=&{\rm ch}\frac{P_0}{\k}+\frac{1}{2k^2}e^{P_0/\k}\vec P^2
\,,\nonumber\\
f^0_{\p m}&=&-\frac{1}{\k}P_m\,,\qquad f^{m}_{\p 0}=-\frac{1}{\k}
e^{P_0/\k}P_m\,,\qquad f^n_{\p m}=\d_m^n\,,\nonumber\\
f^0_{\p4}&=&[{\rm sh}\frac{P_0}{\k}+\frac{1}{2\k^2}e^{P_0/\k}
\vec P^2]\,,\nonumber\\
f^4_{\p 0}&=&[{\rm sh}\frac{P_0}{\k}-\frac{1}{2\k^2}e^{P_0/\k}
\vec P^2]\,,\\
f^m_{\p 4}&=&-\frac{1}{\k}e^{P_0/\k}P_m\,,\qquad f^4_{\p m}=\frac{1}{\k}
P_m\,,\nonumber\\
f^4_{\p 4}&=&{\rm ch}\frac{P_0}{\k}-\frac{1}{2\k^2}e^{P_0/\k}\vec P^2\,.
\nonumber
\end{eqnarray}
we can prove by direct calculations relations (1.15) for $x^{\mu}$ and then
also relations (1.23)

According to (1.18) operators $f^i_{\p j}$ satisfy the
following system of relations ($i,j,k,l=0,1,2,3,4$)
\begin{equation}
f_k^{\p j}f^k_{\p i}=\d_i^j\,,
\end{equation}
where $f_i^{\p j}=g_{ik}g^{jl}f^k_{\p l}$\,.

The additional system
\begin{equation}
f^j_{\p k}f_i^{\p k}=\d_i^j\,,
\end{equation}
is also valid.

According to reltions $\tau^{i^*}=\tau^i$ ($i=0,1,2,3,4$) we may
write for every $a\in{\cal M}_{\k}$
\begin{equation}
a\tau^i=(\tau^ia^*)^*=(f^i_{\p j}(a^*)\tau^j)^*=\tau^j(f^i_{\p j}(a^*)^*=
f^j_{\p k}(f^i_{\p j}(a^*)^*)\tau^k\,,
\end{equation}
and from eq.(1.25) follows that
\begin{equation}
f^i_{\p j}(a^*)^*=f_j^{\p i}(a)\,.
\end{equation}
We shall use this relation in the last section.

\newsection{Dirac operator and quantum Clifford bundle}

Let us define now following elements of ${\cal M}^*_{\k}$:
\begin{equation}
\partial_0=i\k f^4_{\p0}\,,\qquad\partial_m=i\k f^4_{\p m}\,,\qquad
\partial_4=i\k(f^4_{\p4}-1).
\end{equation}

Using these elements we may write the formula for the external derivation
in the compact form ($i=0,1,2,3,4$)
\begin{equation}
da=\partial_i(a)\tau^i
\end{equation}
According to the Leibnitz rule
\begin{equation}
d(ab)=adb+da\,b\,,
\end{equation}
the following system of relations must be satisfied ($i,j=0,1,2,3,4$):
\begin{equation}
\partial_i(ab)=a\partial_i(b)+\partial_j(a)f^j_{\p i}(b)\,,
\end{equation}
or
\begin{equation}
[\partial_i,a]=\partial_j(a)f^j_{\p i}\,.
\end{equation}
Since all $\partial_i\in{\cal M}_{\k}^*$ this is equivalent to
\begin{equation}
\t(\partial_i)=1\otimes\partial_i+\partial_j\otimes f^j_{\p i}\,.
\end{equation}
So to prove (2.2) we must check it for coordinates $x^{\mu}$ and then prove
(2.6) that can be easely done by strict calculations.

According to the general approach \cite{1}, \cite{2} we suppose now that the
Hilbert space ${\cal H}$ is a subspace of ${\Bbb C}^4\otimes{\cal M}_{\k}$.
We take the Dirac operator in the form ($i=0,...4$)
\begin{equation}
D_{\k}=\gamma^{i}\partial_{i}\,,
\end{equation}
where $\gamma^i$ for $i=0,...3$ are the usual Dirac gamma matrices satisfying
the standard relation:
\begin{equation}
\gamma^{\mu}\gamma^{\nu}+\gamma^{\nu}\gamma^{\mu}=2g^{\mu\nu}\,,
\end{equation}
($g^{\mu\nu}={\rm diag}(1,-1,-1,-1)$ is a standard Minkowski space metric)
and $\gamma_4$ is some undefinite matrix which however may be taken in the
form $\g_4=\l I_4$ where $I_4$ is a unit $4\times4$ matrix or
$\gamma_4=\l\gamma_5$ where $\gamma_5=i\gamma_0\gamma_1\gamma_2\gamma_3$.
The choise $\g_4=0$ corresponds to the Dirac operator suggested in \cite{10}.
In this case connection between $D_{\k}$ and $\Box_{\k}$ has the standard form
\begin{equation}
D_{\k}^2=\Box_{\k}
\end{equation}

The direct application of (0.3) gives according to (2.5) the following
expressions for $\tau^i_c$ corresponding to $\tau^i$ elements of quantum 
Clifford bundle.
\begin{equation}
\tau^i=\gamma^jf_j^{\p i}\,,
\end{equation}
Relations (1.21) are fullfilled also for $\tau_c^i$ and the diagram (0.8) is
commutative.

\newsection{The deformed $U(1)$ gauge theory}

In this section we suppose that all up and down indices take values from 0
to 4.

By analogy with classical case we define gauge potentials as elements
of ${\cal A}$ quantum algebra of functions.
Let us introduce the $U(1)$ gauge field by the gauge connection 1-form
\begin{equation}
\omega=iA_k\tau^k
\end{equation}
where $A_{\mu}$ for $\mu=0,..3$ are deformations of usual potential and
$A_4$ may be interpreted as a spin 0 gauge field. Appearance of such scalar
gauge fields in framework of noncommutative geometry was intensively
studied in \cite{1} and \cite{2}.

According to (0.9), (2.7) and connection between $\omega$ and $V$ we can write
the gauge coupled Dirac equation for massless particle in the form
\begin{equation}
\g^k\big_k\psi=0\,,
\end{equation}
where
\begin{equation}
\big_k=\partial_k+igA_jf^j_{\p k}\,.
\end{equation}
($g$ is a gauge charge).

The transformation law (0.13) gives
\begin{equation}
\tilde A_k=UA_jf^j_{\p k}(U^*)-i/gU\partial_k(U^*)
\end{equation}

Defining the curvature form
\begin{equation}
\Omega=d\omega+g\omega\wedge\omega\,,
\end{equation}
we obtain according to (0.13) the following transformation law for it
\begin{equation}
\tilde\Omega=U\Omega U^*\,.
\end{equation}

Defining field strenght tensor by
\begin{equation}
iF_{jk}\tau^j\wedge \tau^k=\Omega\,,
\end{equation}
or according to (1.16)
\begin{equation}
F_{ij}=\partial_i(A_j)-\partial_j(A_i)+iA_k[f^k_{\p i}(A_j)-f^k_{\p j}(A_i)]\,,
\end{equation}
we can using relations (3.6) (1.21) and (1.22) obtain for them the following 
transformation law
\begin{equation}
\tilde F_{ij}=UF_{kl}f^k_{\p i}f^l_{\p j}(U^*)\,.
\end{equation}

We also may obtain the tensor $F_{ij}$ by commuting covariant derivatives
\begin{equation}
[\big_i,\big_j]=igF_{mn}f^m_{\p i}f^n_{\p j}\,.
\end{equation}
It is easy to prove that the following Bianchi identities
\begin{equation}
[\big_i,[\big_j,\big_k]]+[\big_k,[\big_i,\big_j]]+[\big_j,[\big_k,\big_i]]=0\,.
\end{equation}
are satisfied.

Defining deformed covariant derivatives of the strenght tensor as
\begin{equation}
\big_mF^{mk}=\partial_mF^{mk}+ig(A_jf^j_{\p m}(F^{mk})-
F^{mn}f_m^{\p j}f_n^{\p k}(A_j))\,,
\end{equation}
it is easy to obtain for them the following transformation law
\begin{equation}
\tilde \big_mF^{mk}=U\big_mF^{mn}f_n^{\p k}(U^*)\,.
\end{equation}

In the limit $\k\rightarrow\infty$ as it follows from (1.24) (or generally
from (1.21)) we have $f^m_{\p n}\rightarrow\d^m_n$ so that transformation
laws (3.9) and (3.13) may be considered as deformations of standard
formulas. 
Eqs. (3.11) and
\begin{equation}
\big_mF^{mk}=0\,,
\end{equation}
may be interpreted as ${\cal P}_{\k}$-covariant equations of deformed 
electrodynamics in $\k$-Minkowski space. 

It will be interesting to derive (3.14) according to some kind of variational 
principle. We have no recipe to do it. 
However defining elements
\begin{equation}
C=F^{ij}F^*_{ij}\,,\qquad C_+=F_{ij}f^i_{\p k}f^j_{\p l}(F^{kl})\,,\qquad
C_-=f^i_{\p k}f^j_{\p l}(F^*_{ij})F^{kl^*}
\end{equation}
we have
\begin{equation}
\tilde C=UCU^*\,,\qquad\tilde C_{\pm}=UC_{\pm}U^*\,.
\end{equation}
Really for $\tilde C$
\begin{equation}
\tilde C=UF^{kl}f_k^{\p i}f_l^{\p j}(U^*)f_i^{\p u}f_j^{\p v}(U)F_{uv}U^*\,,
\end{equation}

But since $f^i_{\p j}(1)=\d^i_j$ then according to (0.11) and (1.22)
\begin{equation}
f_k^{\p i}f_l^{\p j}(U^*)f_i^{\p u}f_j^{\p v}(U)=f_k^{\p u}f_l^{\p v}
(U^*U)=\d^u_k\d^v_l\,,
\end{equation}
which proves first equation in (3.8). For $C_+$ we have,
\begin{equation}
\tilde C_+=UF_{ij}^*f^i_{\p k}f^j_{\p l}(F^{uv})f^k_{\p c}f^l_{\p d}
f_u^{\p c}f_v^{\p d}(U^*)\,,
\end{equation}
and relation (3.8) for $C_+$ follows now from (1.26) and commutativity of
${\cal M}_{\k}^*$. The proof for $C_-$ follows now from the relation
\begin{equation}
C_-=C_+^*\,,
\end{equation}
which can be easely proved according to (1.28).

Now in order to find a gauge invariant action from $C$ and $C_{\pm}$
we have by analogy with the undeformed case to take integral over
${\cal M}_{\k}$. This means that there exist a subalgebra $L^1({\cal M}_{\k})$
of ${\cal M}_{\k}$ and a positive linear functional $h :\,L^1({\cal M}_{\k})
\rightarrow {\Bbb C}$. It is natural to suppose that $L^1({\cal M}_{\k})$ is
invariant under ${\cal P}_{\k}$-coaction (1.3),
\begin{equation}
\Phi_R(L^1({\cal M}_{\k}))=L^1({\cal M}_{\k})\otimes {\cal P}_{\k}\,.
\end{equation}
Also it is natural to suppose that $h$ is ${\cal P}_{\k}$ invariant so that for
every $a\in L^1({\cal M}_{\k})$
\begin{equation}
(h\otimes id)\circ\t(a)=h(a)1_{{\cal P}_{\k}}\,.
\end{equation}

Let now $U(1)_{{\cal M}_{k}}$ be the group of all $U\in{\cal M}_{\k}$
satisfying (0.11) and additionally preserving $h$, so that for every 
$a\in L^1({\cal M}_{\k})$
\begin{equation}
UaU^*\in L^1({\cal M}_{k})\,,
\end{equation}
and
\begin{equation}
h(UaU^*)=h(a)\,.
\end{equation}

We see now that expressions $h(C)$, $h(C_{\pm})$ and gauge coupled Dirac
equation (3.2) are invariant under the $U(1)_{{\cal M}_{\k}}$ action.

Since now we have not any construction of $h$ and $U(1)_{{\cal M}_{\k}}$
we can not study completely the case of finite ${\k}$.
However we may try to study the case $\k\rightarrow \infty$.
Considering the following Lagrangian density ${\cal L}=-\frac{1}{4}C$
we get from (1.24), (2.1), (3.5), (3.7) in the limit $\k\rightarrow\infty$
\begin{equation}
{\cal L}=-\frac{1}{4}F_{\mu\nu}F^{\mu\nu}+\frac{1}{2}\partial_{\mu}A_4
\partial^{\mu}A_4\,.
\end{equation}
This express the general statement suggested in \cite{1},\cite{2} that
noncommutativity of space-time produces the appearance of spin 0 gauge fields.

\newsection{Conclusions}

In this paper we have defined Dirac operator on $\k$-Minkowski space
according to the A. Connes scheme. In the special case it coincides
with the one suggested in \cite{10}. We also constructed equations for 
deformed electrodynamics on $\k$-Minkovski space. These equations may be 
concidered as covariant under the action of quantum Poincare group. 
We also mentoined the fact of natural appearance of scalar gauge field
in the theory.

\end{document}